\documentclass[a4paper,12pt]{article}
\usepackage{amsmath, amssymb, latexsym, amscd, amsthm,amsfonts,amstext}
\usepackage[mathscr]{eucal}
\usepackage{graphicx}
\usepackage{subfig}

\numberwithin{equation}{section}
\input{tcilatex}
\begin{document}

\title{ On a problem of S.L. Sobolev}
\author{Michael V. Klibanov \\
%EndAName
Department of Mathematics and Statistics \\
\and University of North Carolina at Charlotte \and Charlotte, NC 28223, USA 
\\
%EndAName
mklibanv@uncc.edu }
\date{}
\maketitle

\begin{abstract}
In 1930 Sergey L. Sobolev \cite{Sob1,Sob2} has proposed a construction of
the solution of the Cauchy problem for the hyperbolic equation of the second
order with variable coefficients in 3-d. Although Sobolev did not construct
the fundamental solution, his construction was modified later by Romanov 
\cite{R2,R3} to obtain the fundamental solution. However, these works impose
a restrictive assumption of the regularity of geodesic lines in a large
domain. In addition, it is unclear how to realize those methods numerically.
In this paper a simple construction of a function, which is associated in a
clear way with the fundamental solution of the acoustic equation with the
variable speed in 3-d, is proposed. Conditions on geodesic lines are not
imposed. An important feature of this construction is that it lends itself
to effective computations.
\end{abstract}

% -------------------------------------------------------


\graphicspath{{Figures/}}

%
%
% \documentclass{amsart}
% %%%%%%%%%%%%%%%%%%%%%%%%%%%%%%%%%%%%%%%%%%%%%%%%%%%%%%%%%%%%%%
% \usepackage{amssymb}
% \usepackage[T1]{fontenc}
% \usepackage[latin1]{inputenc}
% \usepackage{graphicx}
% \usepackage{geometry}
% \usepackage{url}
% \usepackage{epsfig}
% % \usepackage{labelfig}
% \usepackage{verbatim}
% % \usepackage{umlaut}
% \usepackage{euscript}
% \usepackage{afterpage}
% \usepackage{graphics}
% \usepackage{amsmath}
% \usepackage{pst-plot}
% \usepackage{subfig}
%
% % \input{psfig.sty}
% \newtheorem{theorem}{Theorem}[section]
% \newtheorem{lemma}[theorem]{Lemma}
% \newtheorem{proposition}[theorem]{Proposition}
% \newtheorem{corollary}[theorem]{Corollary}
% \newtheorem{definition}[theorem]{Definition}
% \numberwithin{equation}{section}
% \newcommand{\norm}[1]{\left\Vert#1\right\Vert}
% \newcommand{\abs}[1]{\left\vert#1\right\vert}
% \newcommand{\set}[1]{\left\{#1\right\}}
% \newcommand{\Real}{\mathbb R}
% \newcommand{\eps}{\varepsilon}
% \newcommand{\To}{\longrightarrow}
% \newcommand{\bn}{\mathbf{n}}
% \def\OFEM{\Omega_{FEM}}
% \def\OFDM{\Omega_{FDM}}
% \newcommand\scal[1]{(\!(#1)\!)}
% \newcommand\bscal[1]{\big(\!\big(#1\big)\!\big)}
% % \input{tcilatex}
% \def\bR{\mathbb{R}}
% \def\bx{\mathbf{x}}
%
%
%
% \begin{document}
%
%
%
%
%

%
% \thispagestyle{plain}
%

\graphicspath{{FIGURES/}
{Figures/}
{FiguresJ/newfigures/}
{pics/}}

\section{Introduction}

\label{sec:1}

In 1930 Sergey L.\ Sobolev, one of the most distinguished mathematicians of
the 20$^{\text{st}}$ century, has published two papers \cite{Sob1,Sob2}
where he has constructed the solution of the hyperbolic equation of the
second order in the 3-d case with variable coefficients in the principle
part of the operator. It was assumed that these coefficient depend on
spatial variables. This result is readily available in the textbook of
Smirnov \cite{Smirnov}. Sobolev did not find the fundamental solution. The
main reason of this was that the notion of the fundamental solution was
unknown in 1930. Still, in his books about inverse problems \cite{R2,R3}
Romanov has modified the method of Sobolev to construct the fundamental
solution of that equation. However, the constructions of both Sobolev and
Romanov impose a quite restrictive assumption on the variable coefficients
in the principal part of the hyperbolic operator. Let the time variable $%
t\in \left( 0,T\right) .$ It is assumed in above cited publications that
geodesic lines generated by these coefficient are regular in a domain $%
Q\left( T\right) \subset \mathbb{R}^{3}.$ The larger $T$ is, the larger $%
Q\left( T\right) $ is. So, $Q\left( \infty \right) =\mathbb{R}^{3}.$ The
regularity of geodesic lines in the domain $Q\left( T\right) $ means that
for any two points $x,y\in Q\left( T\right) $ there exists a single geodesic
line connecting them.

Another construction of the fundamental solution of that equation can be
found in the book of Vainberg \cite{V}. The construction of \cite{V} imposes
the non-trapping condition on variable coefficients in the principal part of
the hyperbolic operator. In addition, the technique of \cite{V} relies on
the canonical Maslov operator, which is not easy to obtain explicitly. Thus,
even though the structure of the fundamental solution in any of above
constructions can be seen, formally at least, still many elements of this
structure cannot be expressed via explicit formulas. For example, neither
the solution of the eikonal equation of the method of Sobolev, nor the
Maslov construction cannot be expressed via explicit formulas.

The numerical factor is important nowadays.\ However, it is not immediately
clear how to compute numerically fundamental solutions obtained in above
references. Therefore, it is also unclear how to compute solutions of Cauchy
problems for heterogeneous hyperbolic equations if using those fundamental
solutions.

In this paper we propose a simple method of the construction of a function,
which is associated in a clear way with the fundamental solution for the
acoustic equation in the 3-d case with the variable coefficient. We impose
almost minimal assumptions on this coefficient. Geodesic lines are not used.
It is important that our function can be both accurately and effectively
approximated numerically via the Galerkin method as well as via the Finite
Difference Method. Thus, we show that this function can be straightforwardly
used for computations of the Cauchy problem for the heterogeneous acoustic
equation.

The author was prompted to work on this paper while he was working on
publications \cite{KR1,KR2} about reconstruction procedures for phaseless
inverse scattering problems. Indeed, in \cite{KR1,KR2} the structure of the
fundamental solution of the acoustic equation in time domain with the
variable coefficient in its principal part is substantially used.

\section{Construction}

Below $x\in \mathbb{R}^{3}$ and $\Omega \subset \mathbb{R}^{3}$ is a bounded
domain. Let $c_{0}$ and $c_{1}$ be two constants such that $0<c_{0}\leq
c_{1}.$ We assume that the function $c\left( x\right) $ satisfies the
following conditions%
\begin{equation}
c\in C^{1}\left( \mathbb{R}^{3}\right) ,c\in \left[ c_{0},c_{1}\right] ,
\label{2.1}
\end{equation}%
\begin{equation}
c\left( x\right) =1\text{ for }x\in \mathbb{R}^{3}\diagdown \Omega .
\label{2.3}
\end{equation}

Consider the following Cauchy problem%
\begin{equation}
c\left( x\right) y_{tt}-\Delta _{x}y=f\left( x,t\right) ,\left( x,t\right)
\in \mathbb{R}^{3}\times \left[ 0,\infty \right) ,  \label{2.15}
\end{equation}%
\begin{equation}
y\mid _{t=0}=y_{t}\mid _{t=0}=0.  \label{2.16}
\end{equation}%
It is assumed here that $f$ is an appropriate function such that there is a
guarantee of the existence of the unique solution $y\in H^{2}\left( \mathbb{R%
}^{3}\times \left( 0,T\right) \right) ,\forall T>0$ of this problem, see,
e.g. Theorem 4.1 in \S 4 of Chapter 4 of the book of Ladyzhenskaya \cite{Lad}
for sufficient conditions for the latter. We will specify $f$ later. The
fundamental solution for the operator $c\left( x\right) \partial
_{t}^{2}-\Delta $ is such a function $P\left( x,\xi ,t,\tau \right) $ that
the solution of the problem (\ref{2.15}), (\ref{2.16}) can be represented in
the form%
\begin{equation}
y\left( x,t\right) =\dint\limits_{0}^{\infty }\dint\limits_{\mathbb{R}%
^{3}}P\left( x,\xi ,t,\tau \right) f\left( \xi ,\tau \right) d\xi d\tau .
\label{2.30}
\end{equation}%
Below we modify formula (\ref{2.30}).

Let $\xi \in \mathbb{R}^{3}$ and $\tau \geq 0$ be parameters. Consider the
following Cauchy problem%
\begin{equation}
c\left( x\right) u_{tt}=\Delta _{x}u+\delta \left( x-\xi \right) \delta
\left( t-\tau \right) ,  \label{2.4}
\end{equation}%
\begin{equation}
u\mid _{t=0}=u_{t}\mid _{t=0}=0.  \label{2.5}
\end{equation}

\subsection{Heuristic part of the construction}

It is convenient for us to work in this subsection with a purely heuristic
derivation. Represent the solution of the problem (\ref{2.4}), (\ref{2.5})
as $u=u_{0}+v$, where $u_{0}$ is the fundamental solution of the wave
equation,%
\begin{equation}
u_{0}=\frac{\delta \left( t-\tau -\left\vert x-\xi \right\vert \right) }{%
4\pi \left\vert x-\xi \right\vert }.  \label{2.6}
\end{equation}%
Thus,%
\begin{equation}
\partial _{t}^{2}u_{0}=\Delta _{x}u_{0}+\delta \left( x-\xi \right) \delta
\left( t-\tau \right) ,  \label{2.60}
\end{equation}%
\begin{equation}
u_{0}\mid _{t=0}=u_{0t}\mid _{t=0}=0.  \label{2.61}
\end{equation}%
Hence, the function $v$ satisfies the following conditions%
\begin{equation}
c\left( x\right) v_{tt}=\Delta _{x}v-\left( c\left( x\right) -1\right) \frac{%
\delta ^{\prime \prime }\left( t-\tau -\left\vert x-\xi \right\vert \right) 
}{4\pi \left\vert x-\xi \right\vert },  \label{2.7}
\end{equation}%
\begin{equation}
v\mid _{t=0}=v_{t}\mid _{t=0}=0.  \label{2.8}
\end{equation}

Consider the operator $A$, 
\begin{equation*}
A\left( f\right) =\dint\limits_{0}^{t}f\left( y\right) dy
\end{equation*}%
for appropriate functions $f$. Purely heuristically again apply the operator 
$A^{4}$ to both sides of equation (\ref{2.7}).\ Denote $\widetilde{w}\left(
x,\xi ,t,\tau \right) =A^{4}\left( v\right) .$ Then (\ref{2.7}) and (\ref%
{2.8}) imply that 
\begin{equation}
c\left( x\right) \widetilde{w}_{tt}-\Delta _{x}\widetilde{w}=-\left( c\left(
x\right) -1\right) \frac{\left( t-\tau -\left\vert x-\xi \right\vert \right) 
}{4\pi \left\vert x-\xi \right\vert }H\left( t-\tau -\left\vert x-\xi
\right\vert \right) ,  \label{2.9}
\end{equation}%
\begin{equation}
\widetilde{w}\mid _{t=0}=\widetilde{w}_{t}\mid _{t=0}=0,  \label{2.10}
\end{equation}%
where $H\left( z\right) $ is the Heavyside function,%
\begin{equation*}
H\left( z\right) =\left\{ 
\begin{array}{c}
1,z>0, \\ 
0,z<0.%
\end{array}%
\right.
\end{equation*}

\subsection{Rigorous part of the construction}

Starting from this point, we are not acting heuristically anymore. To the
contrary, we work rigorously everywhere below.

Consider the Cauchy problem (\ref{2.9}), (\ref{2.10}). Denote%
\begin{equation*}
g\left( x,\xi ,t,\tau \right) =-\left( c\left( x\right) -1\right) \frac{%
\left( t-\tau -\left\vert x-\xi \right\vert \right) }{4\pi \left\vert x-\xi
\right\vert }H\left( t-\tau -\left\vert x-\xi \right\vert \right) .
\end{equation*}%
It follows from (\ref{2.3}) that for any fixed pair $\left( \xi ,\tau
\right) \in \mathbb{R}^{3}\times \left[ 0,\infty \right) $ and for any $T>0$
functions $g,g_{t}\in L_{2}\left( \mathbb{R}^{3}\times \left( 0,T\right)
\right) .$ Hence, Theorem 4.1 of \S 4 of Chapter 4 of the book of
Ladyzhenskaya \cite{Lad}, (\ref{2.1}), (\ref{2.3}) as well as other results
of that chapter imply that for any fixed pair $\left( \xi ,\tau \right) \in 
\mathbb{R}^{3}\times \left[ 0,\infty \right) $ and for any $T>0$ there
exists unique solution $\widetilde{w}\in H^{2}\left( \mathbb{R}^{3}\times
\left( 0,T\right) \right) $ of the Cauchy problem (\ref{2.9}), (\ref{2.10}).
Furthermore, it was shown in the proof of that theorem of \cite{Lad} that
this solution can be effectively constructed numerically via the Galerkin
method. It is also well known that it can be numerically constructed via the
Finite Difference Method. In addition, the energy estimate implies that 
\begin{equation}
\widetilde{w}\left( x,\xi ,t,\tau \right) =0\text{ for }t\leq \tau .
\label{2.11}
\end{equation}%
Also, it follows from (\ref{2.9}) that $\widetilde{w}=\widetilde{w}\left(
x,\xi ,t-\tau \right) $

Consider now the function $w$ defined as%
\begin{equation*}
w\left( x,\xi ,t-\tau \right) =\widetilde{w}\left( x,\xi ,t-\tau \right)
+A^{4}\left( u_{0}\right) .
\end{equation*}%
Hence,%
\begin{equation}
w\left( x,\xi ,t-\tau \right) =\widetilde{w}\left( x,\xi ,t-\tau \right) +%
\frac{\left( t-\tau -\left\vert x-\xi \right\vert \right) ^{3}}{6\cdot 4\pi
\left\vert x-\xi \right\vert }H\left( t-\tau -\left\vert x-\xi \right\vert
\right) .  \label{2.12}
\end{equation}%
Using (\ref{2.6})-(\ref{2.61}), (\ref{2.9}), (\ref{2.10}) and (\ref{2.12})
and applying direct calculations, we obtain that the function $w$ satisfies
the following conditions%
\begin{equation}
c\left( x\right) w_{tt}-\Delta _{x}w=\delta \left( x-\xi \right) \frac{%
\left( t-\tau \right) ^{3}}{6}H\left( t-\tau \right) ,  \label{2.13}
\end{equation}%
\begin{equation}
w\mid _{t=0}=w_{t}\mid _{t=0}=0.  \label{2.14}
\end{equation}

Consider now an arbitrary function $f\left( x,t\right) $ such that 
\begin{equation*}
f\left( x,t\right) \in C^{4}\left( \mathbb{R}^{3}\times \left[ 0,\infty
\right) \right) ,
\end{equation*}%
\begin{equation*}
\partial _{t}^{k}f\left( x,0\right) =0,k=0,1,2,3,4,
\end{equation*}%
\begin{equation*}
f\left( x,t\right) =0,\forall x\in \mathbb{R}^{3}\diagdown G_{f},
\end{equation*}%
where $G_{f}$ is a bounded domain depending on the function $f$. Consider
the function $p_{f}\left( x,t\right) $ defined as%
\begin{equation}
p_{f}\left( x,t\right) =\dint\limits_{0}^{t}\dint\limits_{G_{f}}w\left(
x,\xi ,t-\tau \right) f\left( \xi ,\tau \right) d\xi d\tau .  \label{2.17}
\end{equation}%
Using the integration by parts, (\ref{2.11}) and (\ref{2.12}), we obtain
that there exist four derivatives of the function $p_{f}\left( x,t\right) $
with respect to $t$. In particular, 
\begin{equation}
\partial _{t}^{4}p_{f}\left( x,t\right)
=\dint\limits_{0}^{t}\dint\limits_{G_{f}}w\left( x,\xi ,t-\tau \right)
\partial _{\tau }^{4}f\left( \xi ,\tau \right) d\xi d\tau .  \label{2.18}
\end{equation}%
Next, applying the operator $c\left( x\right) \partial _{t}^{2}-\Delta $ to
both sides of (\ref{2.18}) and using (\ref{2.11}) and (\ref{2.13}), we obtain%
\begin{equation}
\left( c\left( x\right) \partial _{t}^{2}-\Delta \right) \partial
_{t}^{4}p_{f}\left( x,t\right) =\dint\limits_{0}^{t}\frac{\left( t-\tau
\right) ^{3}}{6}\partial _{\tau }^{4}f\left( x,\tau \right) d\tau .
\label{2.19}
\end{equation}%
Using integration by parts in (\ref{2.19}) as well as (\ref{2.141}), we
obtain%
\begin{equation}
\left( c\left( x\right) \partial _{t}^{2}-\Delta \right) \partial
_{t}^{4}p_{f}\left( x,t\right) =f\left( x,t\right) .  \label{2.20}
\end{equation}%
In addition, it follows from (\ref{2.11}), (\ref{2.12}) and (\ref{2.18})
that 
\begin{equation*}
\partial _{t}^{4}p_{f}\mid _{t=0}=\partial _{t}^{5}p_{f}\mid _{t=0}=0.
\end{equation*}%
Since the function $\partial _{t}^{4}p_{f}\in H^{2}\left( \mathbb{R}%
^{3}\times \left( 0,T\right) \right) ,\forall T>0,$ then the uniqueness
theorem for the problem (\ref{2.15}), (\ref{2.16}) and (\ref{2.18}) imply
that 
\begin{equation}
y\left( x,t\right) =\dint\limits_{0}^{t}\dint\limits_{G_{f}}w\left( x,\xi
,t-\tau \right) \partial _{\tau }^{4}f\left( \xi ,\tau \right) d\xi d\tau .
\label{2.200}
\end{equation}%
Thus, the solution of the problem (\ref{2.15}), (\ref{2.16}) is constructed
in (\ref{2.200}) using the function $w$. The formula (\ref{2.200}) is our
first analog of the formula (\ref{2.30}). The integration by parts in (\ref%
{2.200}) and the use of (\ref{2.11}) and (\ref{2.12}) leads to the second
analog of the formula (\ref{2.30}), 
\begin{equation}
y\left( x,t\right) =\dint\limits_{0}^{t}\dint\limits_{G_{f}}\partial
_{t}^{2}w\left( x,\xi ,t-\tau \right) \partial _{\tau }^{2}f\left( \xi ,\tau
\right) d\xi d\tau .  \label{2.26}
\end{equation}

\subsection{Numerical comments}

Recall that the function $w$ can be accurately numerically approximated via
either the Galerkin method or the Finite Difference Method. The derivatives
of the function $f$ can be found analytically, if $f$ is given by an
explicit formula.\ However, if $f$ is given with a noise, then a
regularization method should be applied, see, e.g. the book of Tikhonov and
Arsenin \cite{TA}. In particular, it is explained in this book how to stably
differentiate noisy functions using the regularization.\ Furthermore,
numerical examples of stable computations of first and second derivatives
are presented in \cite{TA}. The form (\ref{2.26}) might be sometimes more
convenient than the form (\ref{2.200}) since second derivatives of noisy
functions are obviously more stable to calculate than fourth derivatives.

In summary, formulas (\ref{2.200}) and (\ref{2.26}) imply that the solution $%
y\left( x,t\right) $ of the problem (\ref{2.15}), (\ref{2.16}) can be
effectively numerically calculated via three steps:

\begin{enumerate}
\item \emph{Step 1}. Compute the function $w$ either via the Galerkin method
or via the Finite Difference method.

\item \emph{Step 2}. Compute corresponding $t-$derivatives of the function $%
f $ should be computed either analytically, if $f$ is given by an explicit
formula, or numerically using the regularization, if $f$ is given with a
noise.

\item \emph{Step 3}. Apply one of formulas (\ref{2.200}) or (\ref{2.26}).
\end{enumerate}

\begin{center}
\textbf{Acknowledgments}
\end{center}

This work was supported by the US Army Research Laboratory and US Army
Research Office grant W911NF-15-1-0233 as well as by the Office of Naval
Research grant N00014-15-1-2330.

\end{document}